\documentclass[8pt,a4paper]{article}
\usepackage{amsfonts}
\usepackage{amsmath}
\usepackage{amssymb}
\usepackage{amsthm}
\usepackage[dvips]{graphicx}
\usepackage{latexsym}
\usepackage{mathrsfs}
\usepackage{amscd}
\usepackage{indentfirst}
\usepackage{textcomp}
\usepackage[numbers,square,sort&compress]{natbib} 
\usepackage[colorlinks=true,allcolors=black ]{hyperref}
\allowdisplaybreaks[4]

\newtheorem{Theorem}{Theorem}[section]

\newtheorem{Lemma}[Theorem]{Lemma}

\theoremstyle{definition}
\newtheorem{Definition}{Definition}[section]
\numberwithin{equation}{section}

\newcommand{\BF}{\mathbb{F}}

\begin{document}
\title{{\bf  Super-biderivations and linear super-commuting maps on simple Lie superalgebras}}
\author{\normalsize \bf Da  Xu$^1$\,\,\,\,Qiyuan  Wang$^1$\,\,\,\,Xiaoning  Xu$^1$}
\date{{{\small{  1. School of Mathematics, Liaoning University, Shenyang, 110036,
China  }}}}
\maketitle

\begin{abstract}
Let $L$ be a finite dimensional simple Lie superalgebra over an algebraically closed field of characteristic different from 2. In this paper, we prove that each skew-supersymmtric super-biderivation of $L$ is inner. Furthermore, we prove that each linear super-commuting map on $L$ is a scalar multiplication map.
\end{abstract}

\renewcommand{\thefootnote}{\fnsymbol{footnote}}
\footnote[0]{ Project Supported by National Natural Science
Foundation of China (No.11501274) and the Science Research Project of Liaoning Provincial Education Department, China (No.L2015203).
\\ Author Email: lnuxxn@163.com (X. Xu)}

\section{Introduction}

In recent years, many researchers were interested in biderivations \cite{Benkovic2009,Bresar1995,Bresar2018,chang20191,chang20192,chen2016,cheng2017,du2013,Ghosseiri2013,han2016,tang2017,tang20181,tang20182,tang20183,wang2013,wang2011,zhang2006}. In \cite{Bresar1995}, Bre\v{s}ar showed that every biderivation of right ideals of prime rings is inner and obtained some conclusions about commuting maps and skew-commuting maps. In \cite{Benkovic2009}, Benkovi\v{c} defined the concept of an extremal biderivation, and proved that under certain conditions a biderivation of a triangular algebra is a sum of an extremal and an inner biderivation. Thereafter, Du \cite{du2013} and Wang \cite{wang2011} obtained similar results for generalized matrix algebras and parabolic subalgebras of simple Lie algebras respectively. Since then biderivations of Lie algebras have been studied widely \cite{han2016,tang2017,tang20181,tang20182,tang20183,wang2013}. Furthermore, in \cite{chen2016,cheng2017} the authors proved that each anti-symmetric biderivation is  inner for simple generalized Witt algebras and the Lie algebra $\mathfrak{gca}$. As an application of biderivations, they showed that every linear commuting map is a scalar multiplication map in their respective papers. Afterwards, Chang \cite{chang20191,chang20192} obtained similar results for  restricted Cartan-type Lie algebras $W(n;\underline{1})$, $S(n;\underline{1})$ and $H(n;\underline{1})$.

Lie superalgebras as a generalization of Lie algebras came from supersymmetry in mathematical physics. Naturally, super-biderivations of Lie superalgebras have aroused the interest of many scholars \cite{bai2023,chang2021,cheng2019,Dilxat2023,fan2017,li2018,tang2020,xia2016,xu2015,yuan2021,yuan2018,zhao2020}. In \cite{xu2015}, the definition of super-biderivations of Lie superalgebras  was introduced for the first time and some properties of super-biderivations on Heisenberg superalgebras were described. In \cite{bai2023,chang2021,cheng2019,fan2017,xia2016,zhao2020}, the authors proved that all skew-supersymmetric super-biderivations are inner for simple modular Lie superealgebras of Witt type and special type, generalized Witt Lie superalgebras, twisted $N=2$ superconformal algebras,  centerless super-Virasoro algebras, super-Virasoro algebras and contact Lie superalgebras. Moreover, utilizing the conclusions of super-biderivations, linear super-commuting maps were studied in \cite{bai2023,fan2017,xia2016}. However, by means of the complete classification of classical simple Lie superalgebras over the complex field $\mathbb{C}$, Yuan \cite{yuan2018} proved that all super-biderivations of classical simple Lie superalgebras without the restriction of skew-supersymmetry are inner. Meanwhile, similar results were obtained by them for Cartan type Lie superalgebras over the complex field $\mathbb{C}$ \cite{yuan2021}. As a generalization of the results of Bre\v{s}ar and Zhao on Lie algebras \cite{Bresar2018}, Tang described the intrinsic connections among linear super-commuting maps, super-biderivations and centroids for Lie superalgebras satisfying certain assumptions \cite{tang2020}. 

This paper is devoted to studying the skew-supersymmtric super-biderivations and linear super-commuting maps of a finite dimensional simple Lie superalgebra, which is denoted by $L$, over an algebraically closed field of characteristic different from 2. And the paper is arranged as follows. In Section 2, we review the basic definitions concerning skew-supersymmtric super-biderivations, super-commuting maps and centroids of a Lie superalgebra. In Section 3, we prove that each skew-supersymmtric super-biderivation of $L$ is inner. In Section 4, we prove that each linear super-commuting map on $L$ is a scalar multiplication map.

\section{Preliminaries}

In this section, we recall some basic notations concerning Lie superalgebras \cite{Scheunert1979,Kac1977}.

Suppose $\mathbb{F}$ is an algebraically closed field of characteristic different from 2. Let  $L=L_{\bar{0}} \oplus L_{\bar{1}}$ be a Lie superalgebra over $\mathbb{F}$. A linear map $f : L \rightarrow L$ is said to be homogeneous of degree $\gamma$, $\gamma\in\mathbb{Z}_2$, if
  \[f(L_\alpha)\subset L_{\alpha+\gamma},\]
for any $\alpha\in\mathbb{Z}_2$. And we write $|f|=\gamma$  for the degree of the homogeneous linear map $f$. If $\gamma=0$, then $f$ is called an even linear map. Throughout this paper, we should mention that once the symbol $|f|$ or $|x|$ appears in an expression, it implies that $f$  is a homogeneous linear map or $x$ is a homogeneous element. A linear map $D : L \rightarrow L$ is called a superderivation of $L$, if
  \[D([x,y])=[D(x),y]+(-1)^{|D||x|}[x,D(y)],\]
for any $x,y\in L$. A bilinear mapping $b : L\times L \rightarrow L$ is said to be homogeneous of degree $\gamma$, $\gamma\in\mathbb{Z}_2$, if
  \[b(L_\alpha,L_\beta)\subset L_{\alpha+\beta+\gamma},\]
for any $\alpha,\beta\in\mathbb{Z}_2$. Similarly, we write $|b|=\gamma$  for the degree of the homogeneous bilinear map $b$. If $\gamma=0$, then $b$ is called an even bilinear map. A bilinear map $\phi : L\times L \rightarrow L$ is called a skew-supersymmtric super-biderivation of $L$, if
  \begin{align*}
  	&\phi(x,y)=-(-1)^{|x||y|}\phi(y,x),\\
  	&\phi(x,[y,z])=[\phi(x,y),z]+(-1)^{(|\phi|+|x|)|y|}[y,\phi(x,z)], 
  \end{align*}
for any $x,y,z\in L$. 	If the bilinear map $\phi_{\lambda}$ : $L\times L\rightarrow L$ is defined by
  \[\phi_{\lambda}( x,y ) =\lambda [ x,y],\]
for any $x,y\in L$, where $\lambda \in \BF$. Obviously $\phi_{\lambda}$ is a skew-supersymmtric super-biderivation of $L$. This class of super-biderivations is called inner super-biderivations. 

\begin{Definition}\cite{tang2020}
	For a Lie superalgebra $L=L_{\bar{0}} \oplus L_{\bar{1}}$, an even linear map $f : L \rightarrow L$ is called a linear super-commuting map if
	  \[[f(x),y]=[x,f(y)],\]
	for any $x,y\in L$.
\end{Definition}

\begin{Definition}\cite{zhang2003}\label{def2.2}
	The centroid of a Lie  superalgebra $L=L_{\bar{0}} \oplus L_{\bar{1}}$ is the associative superalgebra
	  \[\Gamma(L)=\{f\in pl(L)\mid [f,adx]=0,\ \text{for any}\ x\in L\}.\]
	Observe that $\Gamma(L)$ is a associative superalgebra over the field $\mathbb{F}$. Then
	   \[\Gamma(L)=\Gamma(L)_{\bar{0}} \oplus \Gamma(L)_{\bar{1}},\] 
	 where $\Gamma(L)_\alpha=\Gamma(L)\bigcap pl_\alpha(L),\alpha\in\mathbb{Z}_2$.
\end{Definition}

\section{Skew-supersymmtric super-biderivations}

 In this section, we establish several lemmas about the skew-supersymmtric super-biderivations of the simple Lie superalgebra $L$. Moreover, we prove that each skew-supersymmtric super-biderivation of $L$ is inner.
 
 \begin{Lemma}\label{lem3.1}
 	Suppose that L is a finite dimensional simple Lie superalgebra over $\mathbb{F}$. Then 
 	  \[\Gamma(L)_{\bar{0}}=\mathbb{F}\cdot \mathrm{Id},\ \Gamma(L)_{\bar{1}}=0,\]
 	where $\mathrm{Id}$ is the identity map on $L$. And further we have $\Gamma(L)=\mathbb{F}\cdot \mathrm{Id}$.
 	\begin{proof}
 		By the Definition \ref{def2.2}, for any $f\in \Gamma(L)$, we have that
 		  \[f\circ adx=(-1)^{|f||x|}adx\circ f,\]
 		for any $x\in L$. Above equation is equivalent to 
 		  \[f([x,y])=(-1)^{|f||x|}[x,f(y)]=[f(x),y],\]
 		for any $x,y\in L$. 
 		
 		(I) Obviously $\Gamma(L)_{\bar{0}}\supset\mathbb{F}\cdot \mathrm{Id}$, we only need to prove $\Gamma(L)_{\bar{0}}\subset\mathbb{F}\cdot \mathrm{Id}$. Set $0\neq f\in \Gamma(L)_{\bar{0}}$. Due to $\mathbb{F}$ is an algebraically closed field, $f$ has a nonzero eigenvalue $\lambda\in\mathbb{F}$. Let $E_\lambda$ denote the corresponding eigen-subspace. For any $x\in L,y\in E_\lambda$, we hvae
 		\begin{align*}
 			f([x,y])=[x,f(y)]=[x,\lambda y]=\lambda[x,y];\\
 			f([y,x])=[f(y),x]=[\lambda y,x]=\lambda[y,x].
 		\end{align*}
 	    According to the above two equations, it follows that $[L,E_\lambda]\subset E_\lambda$ i.e. $E_\lambda$ is an ideal of $L$. Due to $L$ is simple and $E_\lambda\neq0$, we get that $E_\lambda=L$ i.e. $f=\lambda\cdot\mathrm{Id}\in\mathbb{F}\cdot \mathrm{Id}$.
 	  
 	    (II) Set $f\in \Gamma(L)_{\bar{1}}$. For any $x,y\in L$, we have
 	      \[f([x,y])=(-1)^{|x|}[x,f(y)]=[f(x),y].\]
 	    It is easy to see that the kernel of $f$ is an ideal of $L$, hence either $f= 0$ or else $f$ is bijective. On the one hand,
 	      \[f\circ f([x,y])=(-1)^{|x|}f([x,f(y)])=(-1)^{|x|}[f(x),f(y)].\]
 	    On the other hand,
 	      \[f\circ f([x,y])=f([f(x),y])=(-1)^{\bar{1}+|x|}[f(x),f(y)].\]
 	    Comparing two sides of the above two equations, we have 
 	      \[f\circ f([x,y])=-f\circ f([x,y]),\]
 	    for any $x,y\in L$. Due to $L$ is simple and the characteristic of $\mathbb{F}$ is different from 2, we get 
 	      \[f\circ f=0.\]
 	    Thus we have shown that $f=0$.
 	\end{proof}
 \end{Lemma}

\begin{Lemma}\label{lem3.2}
	Suppose that L is a simple Lie superalgebra over $\mathbb{F}$. Then every skew-supersymmtric superbiderivation $\phi : L\times L \rightarrow L$ can be written as
	  \[\phi(x,y)=f_\phi([x,y]),\]
	where $f_\phi\in \Gamma(L)$.
	\begin{proof}
		Because $L$ is a centerless Lie superalgebra and $L=[L,L]$. Then the conclusion holds from Theorem 3.8 in \cite{tang2020}. 
	\end{proof}
\end{Lemma}

\begin{Theorem}\label{the3.3}
	Suppose that L is a finite dimensional simple Lie superalgebra over $\mathbb{F}$. Then each skew-supersymmtric superbiderivation of $L$ is inner.
	\begin{proof}
		Suppose $\phi : L\times L \rightarrow L$ is a skew-supersymmtric superbiderivation of $L$. By Lemma \ref{lem3.1} and Lemma \ref{lem3.2}, for any $x,y\in L$, there is an element $\lambda\in\mathbb{F}$ such that
		  \[\phi(x,y)=f_\phi([x,y])=\lambda[x,y],\]
		where  $\lambda$ is dependent on $\phi$.
	\end{proof}
\end{Theorem}

\section{Linear super-commuting maps}

In this section, by using the conclusion of skew-supersymmtric superbiderivations we prove that each linear super-commuting map on $L$ is a scalar multiplication map.

\begin{Lemma}\label{lem4.1}
	Suppose that L is a Lie superalgebra over $\mathbb{F}$. Let  $f : L \rightarrow L$ be a linear super-commuting map and 
	  \[\phi_f(x,y)=[f(x),y],\]
	for any $x,y\in L$. Then $\phi_f : L\times L \rightarrow L$ is a skew-supersymmtric super-biderivation.
	\begin{proof}
		It is easy to see that $|\phi_f |=|f|=\bar{0}$. For one thing, we get
		  \begin{align*}
		  	\phi_f (x,y)&=[f(x),y]=[x,f(y)]=-(-1)^{|x||y|}[f(y),x]\\
		  	              &=-(-1)^{|x||y|}\phi_f (y,x),
		  \end{align*}
	    for any $x,y\in L$. For another, we have
	      \begin{align*}
	      	\phi_f (x,[y,z])&=[f(x),[y,z]]\\
	      	                   &=[[f(x),y],z]+(-1)^{(|f|+|x|)|y|}[y,[f(x),z]]\\
	      	                   &=[\phi_f (x,y),z]+(-1)^{(|\phi_f |+|x|)|y|}[y,\phi_f (x,z)].
	      \end{align*}
        for any $x,y,z\in L$.
	\end{proof}
\end{Lemma}

\begin{Theorem}
	Suppose that L is a finite dimensional simple Lie superalgebra over $\mathbb{F}$. Then each linear super-commuting map on $L$ is a scalar multiplication map.
	\begin{proof}
		Let  $f : L \rightarrow L$ be a linear super-commuting map. By Theorem \ref{the3.3} and Lemma \ref{lem4.1}, for any $x,y\in L$, there ia an element $\lambda\in\mathbb{F}$ such that
		  \[[f(x),y]=\phi_f (x,y)=\lambda[x,y],\]
		where  $\lambda$ is dependent on $f$. This implies that $[f(x)-\lambda x, y]=0$ for any $y\in L$. Due to $L$ is simple, it is easy to see that
		  \[f(x)=\lambda x,\]
		for any $x\in L$.
	\end{proof}
\end{Theorem}

%

\end{document}